\documentstyle{amsppt}
\input epsf.tex
\magnification=1200
\pagewidth{6.5truein}
\pageheight{9truein}
\baselineskip=20pt
\NoRunningHeads

\def\figure#1#2#3{\vskip6pt\epsfysize#1truein
		\centerline{\epsfbox{#2}}
		\medskip
		\centerline{ {\smc Figure} #3}\vskip6pt}

\topmatter
\title Essential meridional surfaces for tunnel number one knots 
\endtitle
\author Mario Eudave-Mu\~noz
\endauthor
\address Instituto de Matem\'aticas, UNAM, Circuito Exterior,
Ciudad Universitaria, \break  04510 Mexico D.F., MEXICO
\endaddress
\email eudave\@servidor.unam.mx  and \ mario\@matem.unam.mx \endemail

\abstract We show that for each pair of positive integers $g$ and $n$, there are
tunnel number one knots, whose exteriors contain an essential
meridional surface of genus $g$, and with $2n$ boundary components. We also show
that for each positive integer $n$, there are tunnel number one knots whose
exteriors contain $n$ disjoint, non-parallel, closed incompressible surfaces, 
each of genus $n$.
\endabstract

\keywords Tunnel number one knot, essential surface, meridional surface
 \endkeywords

\subjclass 57M25, 57N10 \endsubjclass

\endtopmatter
\document

\medskip
\head {1. Introduction}\endhead
\medskip

In this paper we 
consider essential surfaces, closed or meridional, properly embedded in the
exteriors of tunnel number one knots. The exterior of a knot $k$ is
denoted by $E(k)=S^3-int\,N(k)$. Recall that a knot $k$ in $S^3$ has
tunnel number one if there exists an arc $\tau$ embedded in
$S^3$ with $k\cap \tau=\partial \tau$, such that $S^3-int\,N(k\cup \tau)$
is a genus 2 handlebody. Such an arc is called an unknotting tunnel for
$k$. Equivalently, a knot $k$ has tunnel number one if there
is an arc $\tau$ properly embedded in $E(k)$, such that $E(k)-int\,N(\tau)$ is a
genus 2-handlebody; in general, the unknotting tunnels we consider are of
this type. Sometimes it is convenient to express a tunnel $\tau^\prime$ for a
knot $k$ as $\tau^\prime=\tau_1\cup \tau_2$, where
$\tau_1$ is a simple closed curve and $\tau_2$ is an arc connecting $\tau_1$ and
$\partial N(k)$; by sliding the tunnel we can pass from one expression to the
other.

A surface $S$ properly embedded in a 3-manifold $M$ is essential if it
is incompressible, $\partial$-incompressible, and non-boundary parallel.
A surface properly embedded in the exterior of a knot $k$ is meridional if  each 
component of $\partial S$ is a meridian of $k$. 
Let $M$ be a compact 3-manifold, and let $S$ be a surface in $M$,
either properly embedded or contained in $\partial M$. Let $k$ be 
a knot in the interior of $M$, intersecting $S$ transversely. 
Let $\hat S=S-int\,N(k)$. The surface $\hat S$ is properly embedded
in $M-int\,N(k)$, and its boundary on $\partial N(k)$, if any, 
consists of meridians of $k$. We say that $\hat S$ is meridionally
compressible in $(M,k)$, if there is an embedded disk $D$ in $M$,
intersecting $k$ at most once,
with $\hat S \cap D=\partial D$, so that $\partial D$ is a nontrivial
curve on $\hat S$, and is not parallel to a component of $\partial \hat S$
lying on $\partial N(k)$. Otherwise $\hat S$ is called meridionally
incompressible. In particular if $\hat S$ is meridionally incompressible in
$(M,k)$, then it is incompressible in $M-k$.

Some results are available on incompressible surfaces in tunnel number
one knot exteriors. Regarding meridional surfaces, 
it is shown in [GR] that the exterior of a tunnel number one knot
does not contain any essential meridional planar surface. Another proof of
this fact is given in [M]. This
says that any tunnel number one knot is indecomposable with respect to tangle sum.
Considering closed surfaces, it is shown in [MS] that there are tunnel
number one knots whose complements contain an essential torus, and such knots
are classified. In [E2] it is proved that for each $g\geq 2$, there
exist infinitely many tunnel number one knots whose complements contain
a closed incompressible surface of genus $g$; such surfaces are also
meridionally incompressible.
 
In this paper we prove the following,

\proclaim{Theorem 3.2} For each pair of integers $g\geq 1$ and $n\geq 1$,
there are tunnel number one knots $K$, such that there is an essential
meridional surface $S$ in the exterior of $K$, of genus $g$, and with $2n$
boundary components. Furthermore, $S$ is meridionally incompressible.
\endproclaim

This gives a positive answer to question 1.8 in [GR].
It follows from [CGLS] that any of the knots of Theorem 3.2 also
contains a closed essential surface of genus $\geq 2$. That surface is
obtained by somehow tubing the meridional surface. However such a surface 
will be meridionally compressible.

Combining the construction of [E2,\S 6] with that of Theorem 3.2, we get the 
following,

\proclaim{Theorem 3.3} For each positive integer $n$, there are tunnel number 
one knots $K$, such that in the exterior of $K$  there are $n$ disjoint,
non-parallel, closed incompressible surfaces, each of 
genus $n$. 
\endproclaim

It follows from the construction that one of the surfaces, say  $S_1$, is 
meridionally compressible while the others are meridionally incompressible. 
It follows also
that the surface $S_1$ is the closest to $\partial E(K)$, that is, $S_1$
and $\partial E(K)$ bound a submanifold $M$ which does not contain any 
of the other
surfaces. It follows from [CGLS,2.4.3] that $S_1$ remains incompressible 
after performing any non-integral Dehn surgery on $K$, and then so does 
any of the other surfaces. This fact, Theorem 3.3 and the observation that
the exterior of a tunnel number one knot is a compact 3-manifold
with Heegaard genus 2, imply the following.

\proclaim{Corollary}  For each positive integer $n$, there are closed, 
irreducible $3$-manifolds $M$, with Heegaard genus $2$, such that in $M$ 
there are $n$ disjoint, non-parallel, closed incompressible surfaces, 
each of genus $n$. 
\endproclaim

This corollary improves one of the results of [Q], where it is shown that for
each $n$, there exist closed irreducible 3-manifolds with Heegaard
genus 2 which contain an incompressible surface of genus $n$.

In Theorem 3.3 the genus of the surfaces grows as much as the number of 
surfaces. This fact is essential, i.e., it is not just 
a consequence of the construction method. It follows from the main Theorem 
of the recent paper [ES], that it is impossible for an irreducible 3-manifold
with Heegaard genus $g$, with or without boundary, to contain an
arbitrarily large number of disjoint and closed incompressible surfaces of bounded
genus. 

The idea of the proof of Theorem 3.2 is the following:
Start with a tunnel number one knot $k$, and unknotting tunnel $\tau$,
and a closed incompressible surface in the complement of $k$
which intersects $\tau$ in two points. We know by [MS] and [E2] that such knots
do exist. Now take an iterate of $k$ and $\tau$, i.e., a knot $k^*$ formed
by the union of two arcs $k^*=k_1\cup k_2$, where $k_1=\tau$ and $k_2$
is an arc lying on $\partial N(k)$. Thus $k^*$ intersects $S$ in two points. It
follows that $k^*$ is a tunnel  number one knot (see Lemma 3.1); an unknotting
tunnel $\tau^*$ for $k^*$ is formed by  the union of $k$ and an arc joining $k$ to
a point in $k_1\cap k_2$. Slide 
$\tau^*$ so that it becomes an arc with endpoints on $k^*$, also
denoted by $\tau^*$.
Now take an iterate of $k^*$ and
$\tau^*$; this is a knot $k^{**}$ with tunnel number one which intersects
$S$ in as many points as desired. If $k^*$ and $k^{**}$ satisfy
certain conditions (Theorem 2.1), the surfaces $S_1=S-int\,N(k^*)$
and $S_2=S-int\,N(k^{**})$ are essential meridional surfaces in the 
exterior of $k^*$ and $k^{**}$, respectively.

Throughout, 3-manifolds and surfaces are assumed to be compact,
connected and orientable. If $X$ is contained
in a 3-manifold $M$, then $N(X)$ denotes a regular neighborhood of $X$ in $M$;
if $X$ is contained in a surface $S$, then $\eta(X)$ denotes a regular
neighborhood of $X$ in $S$. $\Delta(\alpha,\beta)$ denotes the
minimal intersection number of two essential simple closed curves on  a torus
$T$.

I am grateful to E. Sedgwick for suggesting that I prove Theorem 3.3.
I am also grateful to the referee, whose many suggestions greatly improved the
exposition of the paper.

\medskip
\head {2. Construction of essential meridional surfaces}\endhead
\medskip

Let $k$ be a knot in $S^3$, and let 
$\tau^\prime=\tau_1\cup \tau_2$ be an unknotting tunnel for $k$, where
$\tau_1$ is a simple closed curve, and $\tau_2$ is an arc with endpoints
in $\partial N(k)$ and
$\tau_1$. Let $S$ be a closed surface of genus $g$ contained in the exterior
of $k$; then $S$ divides $S^3$ into two parts, denoted by $M_1$ and $M_2$,
where, say, $k$ lies in $M_2$. We say that $S$ is special with
respect to $k$ and $\tau^\prime$ if it satisfies:

\roster
\item $\tau_1$ is disjoint from $S$, and $\tau_2$ intersects
$S$ transversely in one point, so $\tau_1$ lies in $M_1$; 

\item $S$ is essential in $E(k)$.

\endroster

This definition is a variation of the one given in [E2,\S 6].

Note that by [MS], [E1], there exist knots 
with these properties when $g=1$; when $g\geq 2$, the existence of knots
like these follows from [E2,6.1]. 
Note that $M_2\cap N(\tau_2)$ is a cylinder $R\cong D^2\times I$,
so that $R\cap S$ is a disk $D_1\cong D^2\times \{1\}$, and $R\cap N(k)$ 
is a disk $D_0\cong D^2\times \{0\}$. Slide $\tau_1$ over $\tau_2$, 
to get an arc $\tau$ with both endpoints on
$D_0\subset \partial N(k)$, so that $\tau \cap M_2$ consists of two 
straight arcs contained in $R$, i.e., arcs which intersect each disk
$D^2\times \{x\}$ transversely in one point. The surface $S$ and
the arc $\tau$ then intersect in two points. The arc $\tau$ has a
neighborhood $N(\tau)\cong D^2\times I$, so that $N(\tau)\cap M_2\subset R$. 

Let $P$ be a solid torus, $D_0$ a disk contained in $\partial P$,
and $\rho=\{\rho_1,\dots,\rho_n\}$, a collection of arcs properly embedded
in $P$, so that its endpoints lie in $D_0$. We say that this forms a toroidal
tangle with respect to $D_0$, and denote it by $(P,D_0,\rho)$.

Recall that the wrapping number of a knot in a solid torus is defined as the
minimal number of times that the knot intersects any meridional disk of such
solid torus. We define the wrapping number of an arc $\rho_i$ in
$P$ as the  wrapping number of the knot obtained by joining the endpoints of
$\rho_i$ with an arc in $D_0$, and then pushing it into the interior of $P$. This
is well defined.

The tangle $(P,D_0,\rho)$ is good if:
\roster

\item Each arc $\rho_i$ has wrapping number $\geq 1$ in $P$, and there is at 
least one arc $\rho_i$ whose wrapping number in $P$ is $\geq 2$;

\item Each arc $\rho_i$ has no local knots, i.e., if a sphere $S$ intersects
$\rho_i$ in two points, then $S$ bounds a ball $B$ such that 
$B\cap \rho_i$ is an unknotted spanning arc.

\endroster

If the tangle $(P,D_0,\rho)$ is good, then $D_0-\partial\rho$ is incompressible
in $P-\rho$, i.e., there is no disk $D$ properly embedded in $P$, 
disjoint from $\rho$, with $\partial D\subset D_0$, and such that 
$\partial D$ is essential in $D_0-\partial\rho$.

Let $A$ be an annulus in $\partial P$, essential in $P$, so that $D_0\subset A$.
The tangle $(P,D_0,\rho)$ is good with respect to $A$ if:
\roster
\item $(P,D_0,\rho)$ is good.
\item No arc $\rho_i$ is isotopic relative to $\partial \rho_i$ to an arc 
$\lambda$ contained in $A$ (ignoring the other arcs).
\endroster

Let $\hat k$ be a knot contained in the interior of $N(k)\cup N(\tau)$. We 
say that $\hat k$ is specially knotted if:

\roster
\item  $\hat k$ intersects the disk $D_0$ transversely in $2n$ points, so that 
$\hat k\cap N(\tau)$ consists of $n$ straight arcs in $N(\tau)$ 
($\cong D^2\times I$);

\item the toroidal tangle $(N(k),D_0,\rho)$ is good, 
where $\rho=\hat k \cap N(k)$;

\item in the case that $k$ is parallel to a curve lying on $S$, assume also
the following:  let $\gamma$
be the curve on $\partial N(k)$ which cobounds an annulus with a curve on $S$,
and so that $\gamma$ meets $D_0$ in one arc. Let $A=\eta(\gamma\cup D_0)$.
Then $(N(k),D_0,\rho)$ is good with respect to $A$,

\endroster
As $N(\tau)\cap M_2\subset R$, it follows that $\hat k\cap R$ consists
of $2n$ straight arcs. So
$\hat k$ intersects $S$ in $2n$ points. Let $\hat S=S\cap E(\hat k)$.
This is a surface properly embedded in $E(\hat k)$,
whose boundary consists of $2n$ meridians of the knot $\hat k$.

\proclaim{Theorem 2.1} Let $k$ be a knot, $\tau^\prime =\tau_1\cup \tau_2$
an unknotting tunnel for $k$, and $S$ a surface which is special with 
respect to $k$ and  $\tau^\prime$. 
Let $\hat k \subset N(k)\cup N(\tau)$ be a knot which is specially knotted. Then
the surface $\hat S=S\cap E(\hat k)$ is an essential meridional surface 
in the exterior of $\hat k$.
Furthermore, if the surface $S$ is meridionally incompressible
in $(S^3,k)$, then $\hat S$ is meridionally incompressible in $(S^3,\hat k)$. 
If $S$ is 
meridionally compressible, but the wrapping number of some arc $\rho_i$ 
in $N(k)$ is $\geq 3$, where $\rho=N(k)\cap \hat k$, then $\hat S$ is meridionally
incompressible. 
\endproclaim

\demo{Proof} To prove that the surface $\hat S$ is essential in the exterior of
$\hat k$, it suffices to show that it is incompressible, because any two-sided,
connected, incompressible surface in an irreducible
3-manifold with incompressible torus boundary must be
$\partial$-incompressible, unless it is a boundary-parallel annulus, which is not
the case here.

Let $S^\prime = (S\cup \partial R -int\,(D_1))$, so $S^\prime$ is
isotopic to $S$, and let $\tilde S=S^\prime\cap E(\hat k)$; then
$\tilde S$ is a surface isotopic to $\hat S$. Denote by $M_1^\prime$ and
$M_2^\prime$ the complementary regions of $\tilde S$ in $E(\hat k)$, where
$\partial N(k)\cap E(\hat k)$ lies in $M_2^\prime$.
Let $T=\partial N(k)- int\,(D_0)$.
This is a once punctured torus, which is properly embedded in $M_2^\prime$,
i.e., $\tilde S\cap T=\partial T=\partial D_0$. 

Let $D$ be a compression disk for $\tilde S$. Suppose first that
it lies in $M_1^\prime$. As $S^\prime$ is essential in $E(k)$, it follows that
$\partial D$ is a trivial curve on $S^\prime$ which bounds a disk 
$D^\prime \subset S^\prime$, and $D\cup D^\prime$ bounds a 3-ball $B$. As
$\partial D$ is supposed to be essential in $\tilde S$, one arc $\alpha$ 
of $\hat k$ contained in $M_1^\prime$ must in fact be contained in the
3-ball $B$. We may assume that $\alpha = \tau$. Note that $\partial D$ must be
isotopic in $S^\prime-\partial \tau$ to $\partial D_0$. Then the tunnel $\tau$
is contained in a 3-ball, which implies that $k$ is the trivial knot. This is a
contradiction.
 
Suppose then that $D$ lies in $M_2^\prime$. Consider the intersection
between $T$ and $D$. If they do not intersect, then there are
two cases: 
(1) $D$ is contained in $N(k)$. In this case $\partial D$ must lie
on $D_0$, which implies that $\partial D$ is trivial on $\tilde S$,
or that $D_0-\rho$ is compressible in $N(k)-\rho$, which
contradicts the hypothesis.
(2) $D$ is disjoint from $N(k)$. One possibility is that $\partial D$ is
isotopic to $\partial D_0$, but in this case the tunnel $\tau$ is, as above,
contained in a 3-ball which is impossible. Otherwise,  by isotoping
$D$ we may assume that $\partial D$ is contained in $S$, and then $D$ is also a
compression  disk for $S$ disjoint from $N(k)$, which contradicts the hypothesis
that $S$ is incompressible in $E(k)$.

Assume then that $D$ and $T$ have nonempty intersection. 
This intersection consists of a finite number of arcs and simple closed curves. 
Assume also that $D$ has 
been chosen, among all compression disks, to have a minimal
number of intersections with $T$. This implies that any curve or arc
of intersection is essential in $T$, for if one curve (arc) is trivial,
then doing surgery on $D$ with the disk bounded by an innermost curve
(outermost arc) we get a disk with fewer intersections with $T$.

Let $\sigma$ be a simple closed curve of intersection which is
innermost in $D$, so it bounds a disk $D^\prime$ whose interior is disjoint
from $T$. If $D^\prime$ lies in $N(k)$, then $\sigma$ is either a 
meridian of $T$, or it is parallel to $\partial T$, but
in both cases it follows that $D_0-\rho$ is compressible in $N(k)-\rho$.
If the interior of $D^\prime$ is disjoint from $N(k)$, then as $k$ is a 
nontrivial knot, $\sigma$ must be trivial on $T$, which contradicts
the choice of $D$.

Assume then that the intersections between $D$ and $T$ consists only of arcs.
Let $\sigma$ be an outermost arc in $D$ which bounds a disk $E$.
Suppose first that $E\subset N(k)$. Then $\partial E=\sigma\cup \delta$,
where $\delta\subset D_0$. It follows that  $\partial E$ is nontrivial on 
$\partial N(k)$, i.e., it is a meridian of $N(k)$, and then each of the 
$\rho_i$ has wrapping number $\leq 1$ in $N(k)$, which contradicts the
hypothesis. So $E$ cannot be contained in $N(k)$.
Again let $\partial E=\sigma\cup \delta$, where $\sigma$ is contained in 
$T$ and $\delta$ in $\tilde S$. As $\sigma$ is nontrivial in $T$ then $\delta$
is also nontrivial in $\tilde S-D_0$. By isotoping $D$ we can ensure that
$\delta\cap \partial R$  consists of two arcs; let $E^\prime \subset \partial R$ 
be a disk containing these arcs in its boundary. Now $E\cup E^\prime$ is an
annulus with one boundary component on $S$, and the other on
$\partial N(k)$. Here we apply [CGLS,2.4.3], where $M$, $S$, $T$, $r_0$ of
that theorem correspond in our notation to $M_2-int\,N(k)$, $S$, $\partial N(k)$,
and the component of $\partial (E\cup E^\prime)$ lying on $\partial N(k)$,
which we denote also by $r_0$. Clearly $S$ compresses after performing meridional
surgery on $\partial N(k)$. Then part (b) of [CGLS,2.4.3] implies that 
$\Delta (\mu,r_0) \leq 1$, where $\mu$ is a meridian of $\partial N(k)$. 
So either $\mu = r_0$, or $r_0$ goes around $\partial N(k)$ once longitudinally.
The first possibility implies that $S$ is meridionally compressible, and the
second one implies that $k$ is parallel to a curve lying on $S$.
So we are done, unless one of these cases happens. Note that each of
these  possibilities excludes the other, for if $k$ is parallel to a curve on $S$,
and $S$ is meridionally compressible, then either $S$ is compressible
or $S$ is isotopic to $\partial N(k)$.

Suppose first that $S$ is meridionally compressible. Let $\sigma$ be an outermost
arc in $D$, which bounds a disk $E$, so that $\partial E=\sigma\cup \delta$, 
where $\sigma$ is contained in 
$T$ and $\delta$ in $\tilde S$. As above, there is a disk 
$E^\prime \subset \partial R$, such that $E\cup E^\prime$ is an annulus with one
boundary component on $S$, and the other is a meridian of
$\partial N(k)$. Consider all the outermost arcs on
$D$; by the argument given above we can assume that any one of them determines a 
curve on $T$ parallel to $\sigma$. Let $F$ be a region on $D$ adjacent to one
of the outermost arcs, so that all of its intersections with $T$, 
except at most one, are outermost arcs. To find such an $F$, take the collection
of arcs in $D$ which are not outermost arcs, and among these choose
one which is outermost.  $F\subset N(k)$, and then either $\partial F$ 
is trivial on $\partial N(k)$, or $\partial F$ is a meridian of $N(k)$. 
$\partial F$ consists of, say, $2m$ consecutive arcs, 
$\partial F = \sigma_1,\delta_1,\dots,\sigma_m,\delta_m$, where
$\sigma_i\subset T$, and $\delta_i \subset D_0$. Then at least $m-1$ of the
arcs are parallel to $\sigma$, say $\sigma_1,\dots,\sigma_{m-1}$. 
If $\sigma_m$ is not parallel to $\sigma$,
then $\partial F$ would go around $N(k)$ once longitudinally, which is impossible,
for $\partial F$ bounds a disk in $N(k)$.
We conclude that all the arcs $\sigma_i$ are parallel in $T$,
as in Figure 1. 

\midinsert
\figure{2.5}{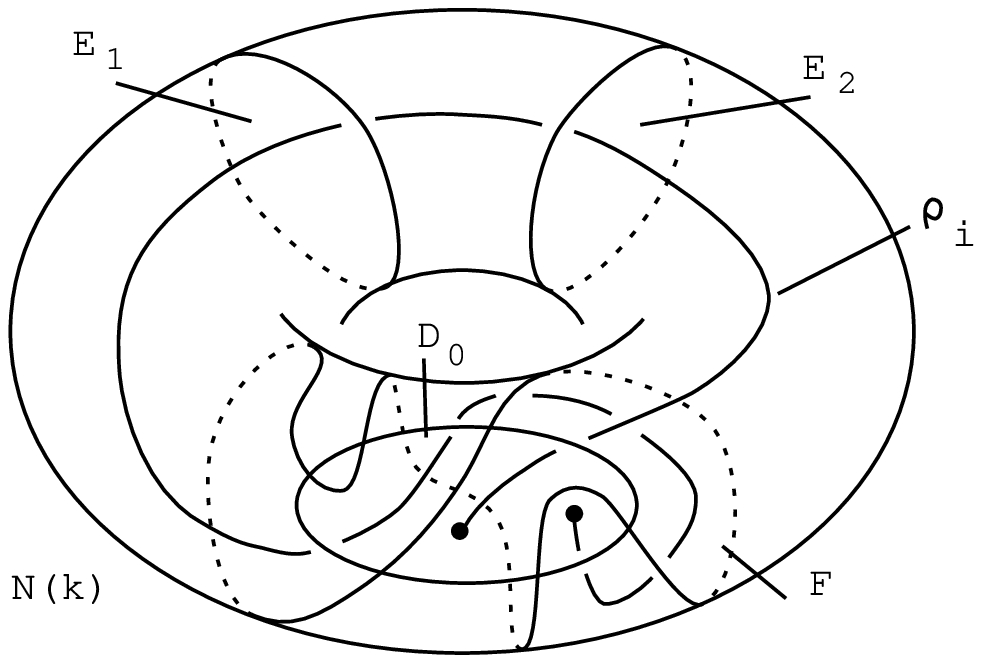}{1}
\endinsert

Let $E_1, E_2$ be two meridian disks of $N(k)$ whose boundaries are disjoint from
$D_0$ and $\cup \sigma_i$. Then $E_1\cup E_2$ bounds a ball $B$ in $N(k)$
which contains $D_0$ and $F$, after possibly isotoping $F$.

There are two cases:

(1) $F$ is parallel to a disk $D_1\subset \partial N(k)$. Clearly
$D_1\subset B$. $F$ and $D_1$ cobound a 3-ball $B_1$. Suppose that an arc
$\rho_i$ is contained in $B_1$. By joining the endpoints of $\rho_i$
with an arc contained in $D_0$, we get a simple closed curve $\rho_i^\prime$,
which is contained in $B$, and then its wrapping number in $N(k)$ is 0,
which contradicts the hypothesis.

So suppose no arc $\rho_i$ is contained in $B_1$. Consider 
$D_1\cap \partial D_0$. This is a collection of arcs which divide $D_1$
into regions which are in $D_0$ or in its complement. If there is
an outermost arc on $D_1$ which bounds a disk contained in $D_0$,
then we can isotope $F$  (and then $D$) through $D_0$ to get a 
compression disk
with fewer intersections with $T$. If no outermost arc  bounds a disk lying in
$D_0$, choose any region $D_0^\prime$  of $D_1\cap D_0$. There is an arc 
$\alpha \subset D_0^\prime$, whose endpoints lie on $\partial D_1$ (then
$\alpha\subset int\,D_0$), and there is a disk $E_0\subset B_1$, so that 
$\partial E_0=\alpha\cup \beta$, where
$\beta$ is an arc on $F$. Cut $D$ along $E_0$, getting two disks;
at least one of them is a compression disk for $\tilde S$, but
it has fewer intersections with $T$.

(2) $\partial F$ is a meridian of $N(k)$, so $\partial F$ is parallel to 
$\partial E_1$ (see Figure 1). So $\partial F$ separates the annulus 
$\partial B-int(E_1\cup E_2)$ into two annuli, denoted by $A_1$ and $A_2$, where
$\partial A_i=\partial E_i\cup \partial F$. Let $\rho_i$ be an arc of
$\rho$, and $\rho_i^\prime$ the simple closed curve
obtained by joining the endpoints of 
$\rho_i$ with an arc in $D_0$. If the endpoints of $\rho_i$ lie in the
same annulus $A_j$, then $\rho_i$ is isotopic rel $\partial \rho_i$ (when
ignoring the other arcs), to an arc disjoint from $E_1$. This implies
that the wrapping number of $\rho_i^\prime$ in $N(k)$ is
0, for $\rho_i^\prime$ is isotopic to a curve disjoint from $E_1$. If the
endpoints of $\rho_i$ lie on different annuli, then $\rho_i$ is isotopic rel
$\partial \rho_i$ to an arc which intersects $E_1$ in one point. This
implies that the wrapping number of $\rho_i^\prime$ in $N(k)$ is 1. This
contradicts the hypothesis that at least one of the arcs have
wrapping number $\geq 2$. This completes the proof when the surface $S$ is
meridionally compressible.

Suppose now that $k$ is parallel to a curve on $S$. As before, let $\sigma$ be an
outermost arc in $D$, which bounds a disk $E$, so that 
$\partial E=\sigma\cup \delta$,  where $\sigma$ is contained in 
$T$ and $\delta$ in $\tilde S$. Recall that the union of $\sigma$ and an arc
on $D_0$ is a curve $\gamma$ on $\partial N(k)$ which cobounds an annulus
$E\cup E^\prime$ with a curve on $S$. Let $A=\eta (\gamma\cup D_0)$. Consider all
the outermost arcs on $D$; recall that any one of them determines a curve on 
$T$ parallel to $\sigma$. Let $F$ be a region on $D$ adjacent to one
of the outermost arcs, so that all of its intersections with $T$, 
except at most one are outermost arcs. 
$F$ is then a disk properly embedded in $N(k)$, which
intersects $D_0$ in $r$ arcs, and all the arcs on $T\cap F$, except at most one
are parallel. 
Let $\partial F=\sigma_1\cup \delta_1 \cup \dots \sigma_r \cup \delta_r$,
where $\sigma_i \subset T$, $\delta_i\subset D_0$, and 
$\sigma_1,\dots \sigma_{r-1}$ are parallel to $\sigma$. There is an annulus
$\Delta$ properly embedded in $N(k)$, $\partial \Delta=\partial A$. We can
assume that $D_0, \sigma_1,\dots, \sigma_{r-1}$ are contained in $A$. If
$\sigma_r$  is not parallel to $\sigma$, then it intersects each component of 
$\partial \Delta$ in one point. It follows that $\partial F$ is trivial in
$\partial N(k)$ if and only if each arc $\sigma_i$ is parallel to $\sigma$.

Suppose first that $\partial F$ is trivial in $\partial N(k)$, then
$\partial F\subset A$, and $F$ is parallel to a disk $D_1\subset A$. We
can assume that $F$ and $\Delta$ do not intersect. $F$ and $D_1$ cobound
a 3-ball $B_1$. Suppose there is an arc $\rho_i\subset B_1$. The arc $\rho_i$
has no local knots, then it is parallel to an arc $\epsilon_i\subset D_1\subset
A$, i.e., the arc $\rho_i$ is isotopic to an arc lying in
$A$, contradicting the hypothesis. See Figure 2.

\midinsert
\figure{2.5}{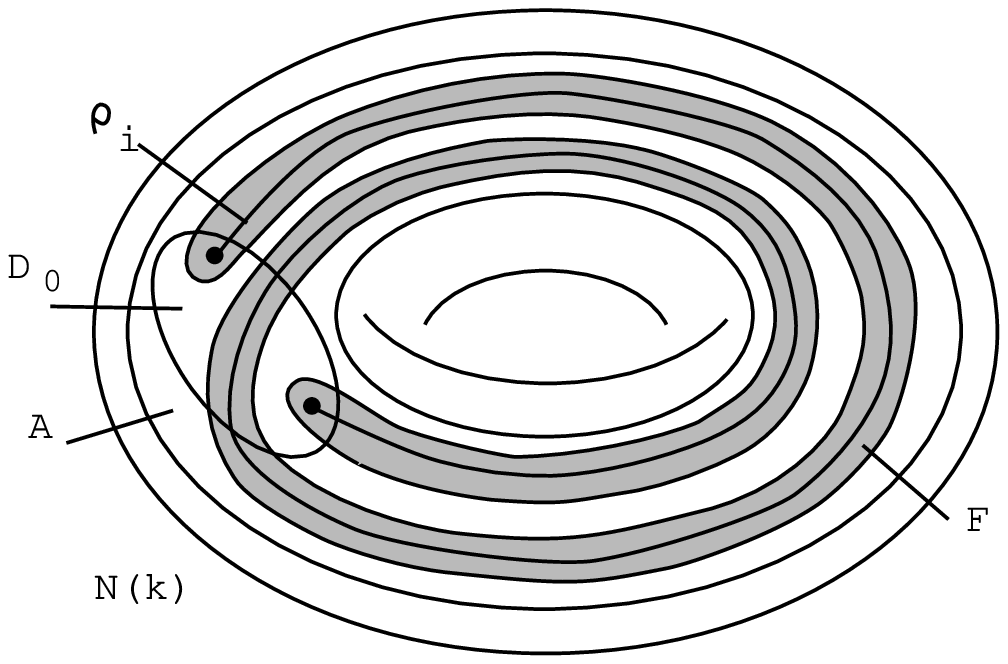}{2}
\endinsert

If there is no arc $\rho_i$ in $B_1$, proceed as in the analogous case when
$\tilde S$ is meridionally compressible, to get a disk $E_0\subset B_1$, with
$\partial E_0= \alpha \cup \beta$, where $\alpha \subset D_0\cap D_1$ and 
$\beta \subset F$, so that by cutting $D$ along $E_0$, we get another compression
disk for $\tilde S$ with fewer intersections with $T$. 

Suppose now that $\partial F$ is a meridian of $N(k)$. Then 
$\partial F=\alpha\cup \beta$, where $\alpha \subset \partial N(k)-A$,
$\beta \subset A$, so $\alpha\subset \sigma_r$.
The annulus $\Delta$ can be isotoped so that $\Delta \cap F$ is a 
single arc. $A$ and $\Delta$ bound a solid torus $\Delta^\prime$, and 
$F\cap \Delta^\prime$ is a meridian disk for $\Delta^\prime$.
If $\rho_i$ is any of the arcs of $\rho$, then $\rho_i$ can be isotoped to
be in the 3-ball $\Delta^\prime-int\,N(F)$, and so it is parallel to an arc
lying on $A$. See Figure 3. This completes the proof.

\midinsert
\figure{2.5}{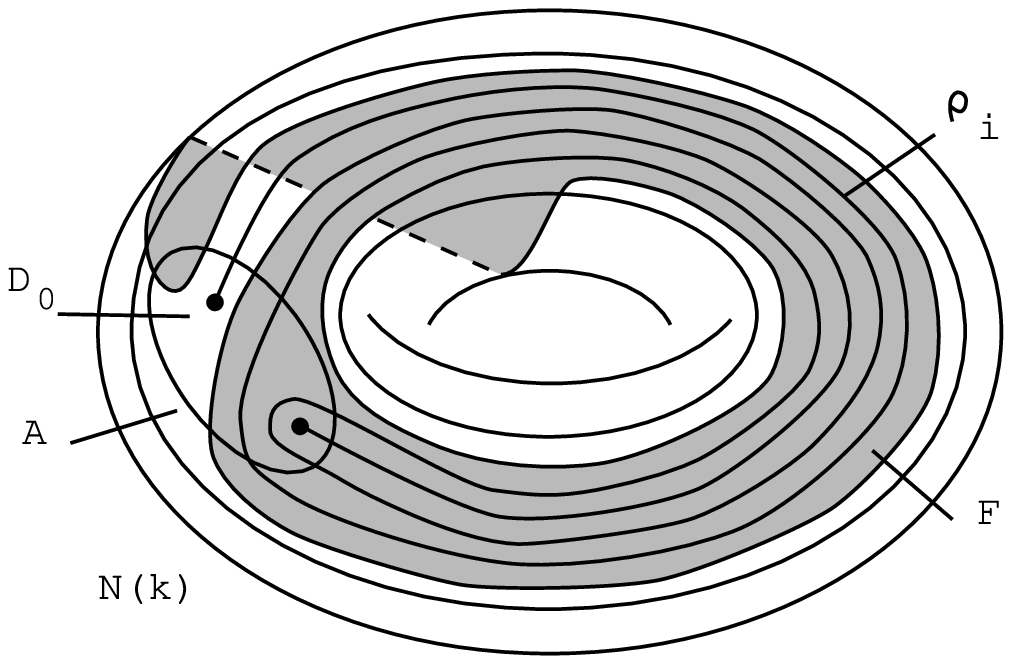}{3}
\endinsert

Now we sketch a proof that 
$\tilde S$ is meridionally incompressible. Suppose there
is a disk $D$ embedded in $S^3$, with $\tilde S\cap D=\partial D$, which
is a nontrivial curve on $\tilde S$, and so that $\hat k$ intersects $D$
transversely in one point. 
If the disk $D$ lies in $M_1^\prime$, then as $S^\prime$ is incompressible in
$E(k)$, it follows that $\partial D$ is a trivial curve on $S^\prime$, which is
the boundary of a disk contained in $S^\prime$ which intersects $\hat k$ once, so
$D$ is not a disk of meridional compression.  

So assume that $D$ lies in $M_2^\prime$.
Look at the intersections between $D$ and
$T$, and suppose $D$ has been chosen to have minimal intersection
with $T$. This implies that any curve or arc of intersection is essential
in $T$. 

Suppose there is a curve of intersection, innermost in
$D$, which bounds a disk $D^\prime$ which meets $\hat k$ once. Then $D^\prime$
must lie in $N(k)$. If $\partial D^\prime$ is a meridian of $N(k)$, then each
$\rho_i$ has wrapping number $\leq 1$. If $\partial D^\prime$ is not a meridian
of $N(k)$, then $\partial D^\prime$ bounds a disk in $\partial N(k)$ which either
lies in $T$ or contains $D_0$. In either case it is impossible for $\rho$ to meet
$D^\prime$ in exactly one point. This shows that simple closed curves of
intersection cannot bound disks which intersect $\hat k$, and then these
curves can be removed as before. 
Suppose there is an outermost arc $\sigma$ in $D$ which bounds a disk $E$ 
disjoint from $\hat k$. Doing an argument as the one done to prove the 
incompressibility of $\tilde S$, we have that $E$ does not lie in $N(k)$.
By the same argument, such a disk can exist only if $S$ is meridionally
compressible, or if $k$ is parallel to a curve on $S$.
Note that it is always  possible to find an outermost arc which bounds a
disk disjoint from $\hat k$. So the proof is complete, except if we have one of
the cases just mentioned.

Suppose first $S$ is meridionally compressible. In this case we suppose that
the wrapping number of some arc $\rho_i$ in $N(k)$ is $\geq 3$. Take an outermost
arc of intersection in $D$, and suppose it bounds a disk
$D^\prime$ contained in $N(k)$, which intersects $\hat k$ in at most one point. 
$\partial D^\prime$ is a meridian of
$N(k)$, and then the wrapping number of any arc $\rho_i$ in $N(k)$ is $\leq 2$,
contradicting the hypothesis in this case. So suppose all outermost
arcs bound disks which do not lie in $N(k)$. As in the proof of the 
incompressibility of $\tilde S$, these arcs in $T$ are all parallel, and each one
of them, together with an arc in $D_0$ is a meridional curve on $\partial N(k)$.
If there is a region
$F\subset D$, such that all the intersections of
$F$ with $T$, except at most one are outermost arcs, and so that $F$
is disjoint from $\hat k$, proceed as in the proof of the incompressibility of
$\tilde S$. If there is no such region $F$, then $T\cap D$ consists of $m$ arcs,
all of which  are outermost arcs in $D$, so that the complement of the arcs
is a single region $F^\prime$ contained in $N(k)$ and intersecting $\hat k$
once. There are two cases:

(1) $F^\prime$ is parallel to a disk $D_1\subset \partial N(k)$. $F^\prime$ and
$D_1$ cobound a 3-ball $B_1$. If some arc $\rho_j$ is contained in $B_1$ then its
wrapping number in $N(k)$ is $0$, contradicting the hypothesis. So there is
just one arc $\rho_i$ which intersects $B_1$; one of its endpoints is in $D_1$
and the arc intersects $F^\prime$ in one point. As $\rho_i$ has no local knots,
$B_1\cap \rho_i$ is an unknotted spanning arc in $B_1$. As in the case of the
incompressibility, there is a disk $E_0\subset B_1$, so that
$\partial E_0=\alpha \cup \beta$, where $\beta$ is an arc on $F^\prime$, and
$\alpha$ is an arc in $D_0 \cap D_1$. Cut $D$ along $E_0$, getting two disks;
at least one of them is a meridional compression disk for $\tilde S$, but
it has fewer intersections with $T$.

(2) $\partial F^\prime$ is a meridian of $N(k)$. The same proof as in the case of
the incompressibility show that if this happens then the wrapping number of any
arc $\rho_i$ in $N(k)$ is $\leq 2$.

Suppose now that $k$ is parallel to a curve on $S$. If $\sigma$ is an outermost
arc of intersection in $D$, bounding a disk $E$ which does not intersect
$\hat k$, then $E$ is not contained in $N(k)$, and 
$\partial E = \sigma\cup \delta$, where 
$\delta$ is contained in $\tilde S$. The union of $\sigma$ and an arc on $D_0$ is
a curve $\gamma$ on $\partial N(k)$ which cobounds an annulus with a curve on $S$;
so $\gamma$ is a curve which goes around $N(k)$ once longitudinally.
If there is an outermost arc of intersection in $D$ bounding a disk $D^\prime$
which intersects $\hat k$ in one point, then $D^\prime$ is a meridian disk of
$N(k)$. In particular this shows that $D\cap T$ cannot consist of just one arc.
As before, if $\sigma^\prime$ is another outermost arc in $D$ which bounds a disk
disjoint from $\hat k$, then $\sigma^\prime$ is parallel to $\sigma$.

Now proceed as in the proof of the incompressibility of
$\tilde S$ in the case that
$k$ is parallel to a curve on $S$. The point is to find a region $F\subset D$,
such that all the intersections of $F$ with $T$, except at most one, are outermost
arcs, and so that $F$ is disjoint from $\hat k$. If such region exists we are
done. If there is an outermost arc on $D$ which bounds a disk which intersects
$\hat k$, then such region $F$ does exists, for otherwise $D\cap T$ will
consist of just one arc. If such a region $F$ does not exist, the only possibility
left is that $D\cap T$ consists of
$m$ arcs, all of which are outermost arcs, and $D\cap N(k)$ is a single disk
$F^\prime$ which meets $\hat k$ once. Then $\partial F^\prime$ is completely
contained in the  annulus $A=\eta(\gamma\cup D_0)$, which implies that 
$\partial F^\prime$ is trivial  on $A$. So $F^\prime$ bounds a disk $D_1$
contained in $A$; $F^\prime$ and $D_1$ cobound a ball $B_1$. If some arc $\rho_j$ 
is contained in $B_1$ then it is parallel to an arc lying on $A$,
contradicting the hypothesis. So there is
just one arc $\rho_i$ which intersects $B_1$; one of its endpoints is in $D_1$
and the arc intersects $F^\prime$ in one point. As $\rho_i$ has no local knots,
$B_1\cap \rho_i$ is an unknotted spanning arc in $B_1$. As in the proof of the
incompressibility of $\tilde S$, we can boundary compress $D$, getting another
meridional compression disk for $\tilde S$, but
with fewer intersections with $T$.
\qed 
\enddemo

\remark{Remark} The conditions imposed on the tangle $(N(k),D_0,\rho)$
are somehow local, i.e., they consider each arc separately. Giving to the tangle
some global property might produce a slightly stronger theorem.
\endremark

\medskip
\subhead {3. Tunnel number one knots and meridional surfaces}\endsubhead
\medskip
Let $k$ be a tunnel number one knot, and $\tau$ an unknotting tunnel for
$k$ which is an embedded arc with endpoints lying on $\partial N(k)$.
Assume that a neighborhood $N(k\cup \tau)$ is decomposed as
$N(k\cup \tau) = N(k) \cup N(\tau)$, where $N(k)$ is a solid torus, 
$N(\tau)\cong D^2\times I$, $N(k)\cap N(\tau)$ consists of two
disks $E_0$ and $E_1$, and $\tau=\{0\}\times I$.

Let $k^*$ be a knot formed by the union of two arcs, 
$k^*=k_1\cup k_2$, such that
$k_1$ is contained in $\partial N(k)$, and $k_2=\tau$.
We say that $k^*$ is an iterate of $k$ and $\tau$. 

\proclaim {Lemma 3.1} Let $k$ and $\tau$ be as above, and let $k^*$ be an iterate
of $k$ and $\tau$. Then $k^*$ is a tunnel number one knot. An unknotting
tunnel $\beta^\prime$ for $k^*$ is given by the union of $k$ and a straight arc
in $N(k)$ connecting $k^*$ and $k$.
\endproclaim

\demo{Proof} $N(k)-k$ is homeomorphic to a product $T\times [0,1)$.
Let $\delta$ be a straight arc in $N(k)$ connecting $k$ and one of the
points $k_1\cap k_2$, i.e., it is an arc which intersects each torus 
$T\times \{x\}$ in one point.
Then $\beta^\prime = k\cup\delta$ is an unknotting tunnel for $k^*$. 
To see that, slide
$k_1$ over $\delta$ and then over $k$ to get a 1-complex which is clearly
equivalent to $k\cup \tau$, so its complement is a genus 2 handlebody.
\qed
\enddemo

Let $k^*$ be an iterate of $k$ and $\tau$. It follows by construction
that $k^*\subset N(k\cup \tau)$. Also if $\beta^\prime$ is the unknotting tunnel
for $k^*$ given by the lemma, then $k^*\cup \beta^\prime \subset N(k\cup \tau)$.
Now $\beta^\prime$ can be modified to be an arc $\beta$ with endpoints in
$k^*$. It follows that if $k^{**}$ is an iterate of $k^*$ and $\beta$, then
$k^{**}$ can be isotoped to lie in $N(k\cup \tau)$. By isotoping $k^{**}$,
if necessary, we have that $k^{**}\cap N(\tau)$ consists of a collection of
arcs parallel  to $\tau$.

\proclaim{Theorem 3.2} For each pair of integers $g\geq 1$ and $n\geq 1$,
there are tunnel number one knots $K$ such that there is an essential
meridional surface $\hat S$ in the exterior of $K$, of genus $g$, and with $2n$
boundary components. Furthermore, $\hat S$ is  meridionally incompressible.
\endproclaim

\demo{Proof}
Let $k$ be a tunnel number one knot. Suppose that $k$ has an unknotting tunnel
$\tau^\prime=\tau_1 \cup \tau_2$, where $\tau_1$ is a simple closed
curve, and $\tau_2$ is an arc connecting $k$ and $\tau_1$. Suppose there is a
closed surface $S$ of genus $g\geq 1$ embedded in the exterior of $k$, which is
special with respect to $k$ and $\tau^\prime$. 

$S$ divides $S^3$ into two parts, $M_1$ and $M_2$, where, say, $\tau_1$
is contained in $M_1$. $M_2\cap N(\tau^\prime)$ is a cylinder 
$R\cong D^2\times I$, so that $R\cap S$ is a disk $D_1$, and $R\cap N(k)$ 
is a disk $D_0$.
Slide $\tau_1$ over $\tau_2$, to get an arc $\tau$ with both endpoints on
$D_0\subset N(k)$, so that $\tau \cap M_2$ consists of two straight arcs
contained in $R$. The surface $S$ and the arc $\tau$ intersect in two points.

Let $k^*$ be an iterate of $k$ and 
$\tau$; then $k^*=k_1\cup k_2$, where $k_2$ is an arc parallel to $\tau$,
so it intersects $S$ in two points. Now $k_1$ is an arc in $\partial N(k)$
whose endpoints lie on $D_0$. By pushing $k_1$ into the interior of $N(k)$
we get a properly embedded arc in $N(k)$. Clearly $k_1$ can be
chosen so that $(N(k),D_0,k_1)$ forms a good tangle, just by taking an arc
whose wrapping number in $N(k)$ is $\geq 2$. Note that $k_1$ has no local knots
in $N(k)$, for it is parallel to an arc lying in $\partial N(k)$. If
$k$ is parallel to a curve on
$S$, let $\lambda$ be the curve on $\partial N(k)$ which cobounds
an annulus with a curve on $S$, so that $\lambda$ meets $D_0$ in one arc;
let $A=\eta(\gamma\cup D_0)$. Clearly $k_1$ can be chosen so that 
$(N(k),D_0,k_1)$ is good with respect to $A$, say by twisting $k_1$
meridionally as many times as necessary; this can be done because the annulus
$A$ goes longitudinally once around $N(k)$, and the wrapping number of $k_1$ in
$N(k)$ is $\geq 2$.
So $k^*$ can be chosen to be specially knotted in $N(k\cup \tau)$. It
follows from Theorem 2.1 that $\hat S=S\cap E(k^*)$ is an essential meridional
surface in
$E(k^*)$, and $\partial \hat S$ consists of two meridians of $k^*$.

This implies that $k^*$ is a tunnel number one knot which has an unknotting
tunnel
$\beta^\prime=k\cup \delta$, where $\delta$ is a straight arc connecting $k^*$
and $k$.  Let $\beta$ be the arc obtained after sliding $k$ over $\delta$.
$N(k^* \cup \beta)$ can be chosen so that it is contained in $N(k\cup \tau)$.
Let $k^{**}$ be an iterate of $k^*$ and
$\beta$, so $k^{**}=\kappa_1\cup \kappa_2$, where $\kappa_1\subset \partial
N(k^*)$, and $\kappa_2$ is the tunnel $\beta$. Note that
$k^{**}\subset N(k\cup \tau)$. The arc $\kappa_1$
can be isotoped so that $\kappa_1\cap N(\tau)$ consist of straight arcs,
and it can be chosen so that $\kappa_1\cap N(\tau)$ consists
of $n$ arcs, $n$ being a fixed positive integer. So $k^{**}$
intersects $S$ in $2n$ points. $k^{**}\cap N(k)$ then consists
of $n$ arcs, $\rho_0,\dots,\rho_n$, which are properly embedded on $N(k)$, 
and whose endpoints lie in $D_0$. Clearly $k^{**}$ can be chosen so that
$(N(k),D_0,\rho)$ forms a good tangle, say by choosing them so that each arc
$\rho_i$, except one, is parallel to $k_1$, and so that each has wrapping number
$\geq 2$. The remaining arc can be chosen to be a band sum of the arc $k_1$ and
the knot
$k$, so it can be chosen to have wrapping number $\geq 3$.
If $k$ is parallel to
$S$, again $k^{**}$ can be chosen so that $(N(k),D_0,\rho)$ is good with
respect to $A$. Then by Theorem 2.1, the surface  $\hat S=S\cap E(k^{**})$ is an
essential meridional surface in $E(k^{**})$, and $\partial \hat S$ consists of
$2n$ meridians of $k^{**}$. 

If $S$ is meridionally incompressible, then $\hat S$ is meridionally
incompressible. If $S$ is meridionally compressible, then $k^*$ and
$k^{**}$ can be chosen so that $\hat S$ is meridionally incompressible.
\qed
\enddemo

\midinsert
\figure{3.3}{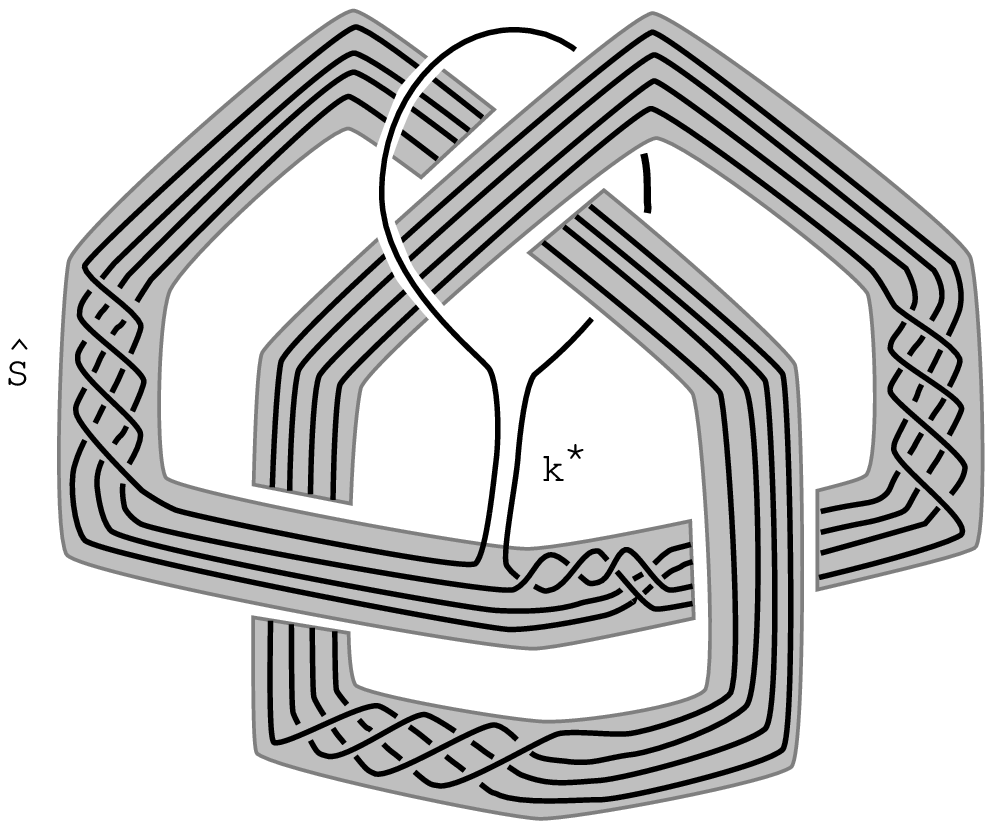}{4}
\endinsert

It follows from the proof of Theorem 3.2 that for the knots $k$ 
constructed in [E2], there are many iterates of $k$,
whose exteriors contain an essential meridional surface. This is because for 
such knots, there is an unknotting tunnel $\tau^\prime$ and a surface $S$
which is special with respect to $k$ and $\tau^\prime$. Note also that some of 
these knots
$k$ are parallel to the surface $S$, while others are not [E2,8.2]. 

An example which illustrates Theorem 3.2 is shown in Figure 4.  Let
$k$ be the (2,-11)-cable of the left hand trefoil; there is a torus
$S$ and unknotting tunnel $\tau^\prime$ for $k$, so that $S$ is special
with respect to $k$ and $\tau^\prime$. Note that $k$ is parallel to a curve on
$S$. The knot $k^*$ shown in Figure 4 is an iterate of $k$ and the tunnel 
$\tau^\prime$. It is not difficult to check that $k^*$ satisfies the conditions
of Theorem 2.1. So it follows that $\hat S$ is an essential meridional surface in
$E(k^*)$.

Combining the last theorem and the construction given in [E2,\S 6],
we get the following.

\proclaim{Theorem 3.3} For each positive integer $n$, there are tunnel number 
one knots $K$, such that in the exterior of $K$  there are $n$ disjoint,
non-parallel, closed incompressible surfaces. Each of the surfaces has
genus $n$. One of the surfaces is meridionally compressible; the
others are  meridionally incompressible.
\endproclaim

\demo{Proof}
Recall the construction given in [E2,\S 6]. Let $k_1$ be a knot,
$\tau^\prime=\tau_1\cup \tau_2$ an unknotting tunnel, and $S_1$ an essential
surface of genus $g$ embedded in $E(k_1)$, which intersects $\tau^\prime$
in one point. So $S_1$ is special with respect to $k_1$ and $\tau^\prime$ (in both
definitions, the one given here and the one in [E2,\S 6]; see [E2,6.1], which
shows that this is true). 
Let $T=\partial N(k_1)$. Let $A$ be an annulus contained in
$T$, and let $\alpha$ be the core of this annulus. Suppose that $\alpha$ wraps 
around $N(k_1)$ at least twice longitudinally. If $k_1$ is parallel to $S_1$, 
suppose also that $\Delta(\gamma,\alpha)\geq 2$,
where $\gamma$ is a curve on $\partial N(k_1)$ which cobounds an annulus
with a curve on $S_1$.

$S_1$ divides $S^3$ into two parts, $M_1$ and $M_2$, where, say, $k_1$
is contained in $M_2$. 
Let $\tau_2^\prime=M_2\cap \tau_2$; so $\tau_2^\prime$ is an arc with an endpoint
on $S_1$ and the other on $\partial N(k)$, which we assume lies on the curve
$\alpha$. The curve $\alpha$ goes around $N(k)$ at least twice longitudinally,
then it is a toroidal graph of type 1 in $N(k)$, as defined in [E2,\S4].
Let $M=M_2-int\,N(k)$. $M$ is a 3-manifold with incompressible boundary.
To show that $\tau_2^\prime\cup \alpha$ is a cabled graph in $M_2$, as defined in
[E2,\S 6], it suffices to prove that $S_1$ remains incompressible after Dehn
filling $M$ along $\partial N(k)$ with slope $\alpha$. If $k$ is not parallel
to a curve on $S_1$, then as $\Delta(\alpha,\mu)\geq 2$, this follows from the
main Theorem of [Wu]. If $k$ is parallel to a curve on $S_1$, then by hypothesis,
$\Delta(\alpha,\gamma)\geq 2$, and by [CGLS,2.4.3] it follows that $S_1$ remains
incompressible.

$N(\tau_2^\prime)$ is a cylinder 
$R\cong D^2\times I$, so that $R\cap S_1$ is a disk $D_1$, and $R\cap N(k_1)$ 
is a disk $D_0$.
Assume that $D_0\subset A$. Consider the manifold
$W=M_1\cup R \cup N(A)$, and let 
$\Sigma=\partial W$. As $\tau_2^\prime\cup \alpha$ is a cabled graph in $M_2$, 
it follows from [E2,6.3] that $\Sigma$ is incompressible in $S^3-int\,W$.

Let $\tau$ be the arc obtained by sliding $\tau_1$ over $\tau_2$,
so that $M_2\cap \tau \subset R$.
Now take an iterate $k_2$ of $k_1$ and $\tau$ of a special form. As before 
$k_2=\kappa_1\cup \kappa_2$, where $\kappa_2=\tau$, and $\kappa_1$ is an 
arc in $\partial N(k_1)$.
Suppose that $\kappa_1$ is contained in $A$, so that its wrapping number 
in $N(A)$ is $\geq 2$ (i.e., $\rho= k_2\cap N(A)$ is a properly embedded arc in
$N(A)$ whose endpoints lie on $D_0^\prime=R\cap \partial N(A)$, and we are
requiring that the curve obtained from $\rho$ by joining its
endpoints with an arc lying on $D_0^\prime$ has wrapping number $\geq 2$ in
$N(A)$). Then
$k_2 \subset W$, and it  follows from [E2,6.4] that $\Sigma$ is incompressible and
meridionally  incompressible in $(W,k_2)$. So $\Sigma$ is a meridionally 
incompressible surface contained in the exterior of $k_2$ of genus $g+1$. By
[E2,8.2], it follows that
$k_2$ is not parallel to a curve lying on $\Sigma$

It is not difficult to see that the knot $k_2$ 
also satisfies the hypothesis of Theorem 2.1; in particular, note that the
arc $\kappa_1$ has wrapping number $\geq 4$ in $N(k_1)$ (for $\kappa_1$ has
wrapping number $\geq 2$ in $N(A)$, and $\alpha$ has winding number $\geq 2$ in
$N(k_1)$). Therefore the surface $\hat S = S_1\cap E(k_2)$
is meridionally incompressible in $E(k_2)$, its boundary consists
of two meridians of $k_2$. By tubing $\hat S$, we
get two  closed surfaces in $E(k_2)$, of
genus $g+1$. By an application of the handle addition Lemma [J], one of the
surfaces must be incompressible in $E(k_2)$; this has to be the surface lying on
$M_2$, for the one lying in $M_1$ bounds a handlebody. Denote by
$\bar S$ such an incompressible surface; note that it is meridionally
compressible.  Then there are two different closed incompressible
surfaces in $E(k_2)$, $\Sigma$ and $\bar S$. By  isotoping $\bar S$ into $W$,
these surfaces become disjoint and are obviously non-parallel.

Note that there is an unknotting tunnel $\beta^\prime=k_1\cup \delta$ for $k_2$,
where $\delta$ is a straight arc in $N(k_1)$ connecting $k_1$ and $k_2$
which intersects both surfaces $\Sigma$ and $\bar S$ in one point. Then
$\Sigma$ and $\bar S$ are both special with respect to $k_2$ and $\beta$. Note
that $\bar S$ is closer to $k_2$ and $\Sigma$ is closer to $k_1$; that is,
the arc $\delta$, when going from $k_2$ to $k_1$, intersects first
$\bar S$ and then $\Sigma$.

We have proved that there is a tunnel number one knot $k_2$ which has an
unknotting tunnel $\tau^\prime=\tau_1\cup \tau_2$, and two disjoint, non-parallel
closed incompressible surfaces in its exterior, each of genus $g+1$, denoted by
$\Sigma$ and $\bar S$, and which are special with respect to $k_2$ and
$\tau^\prime$.
$\Sigma$ is meridionally incompressible and $\bar S$ is meridionally compressible,
and the arc $\tau_2$, when going from $k_2$ to $\tau_1$ intersects first
$\bar S$ and then $\Sigma$. Furthermore, $k_2$ is not parallel to a curve lying 
on any of the two surfaces.

Suppose by induction that we have a tunnel number one knot $k_n$, which has an
unknotting tunnel $\tau^\prime=\tau_1\cup \tau_2$, and $n$ disjoint, non-parallel
closed incompressible surfaces in its exterior, of genus $g+n$, denoted by 
$S_1,S_2,\dots,S_n$, which are special with respect to $k_n$ and $\tau^\prime$.
$S_2,\dots,S_n$ are meridionally incompressible and $S_1$ is
meridionally compressible, and the arc $\tau_2$, when going from $k_n$ to
$\tau_1$ intersects the surfaces in the order $S_1,S_2,\dots,S_n$.
Furthermore, $k_n$ is not parallel to a curve lying  on any
of the  surfaces.

The above construction can be repeated with $k_n$, 
$\tau^\prime=\tau_1\cup \tau_2$ and
$S_1,S_2,\dots,S_n$. $S_i$ divides $S^3$ into $M_1^i$ and $M_2^i$, where
$k_n$ lies in $M_2^i$. Clearly, if $i<j$ then $M_2^i \subset M_2^j$.
Let $\alpha$ be a simple closed curve on $\partial N(k_n)$, which goes at least
twice longitudinally around $N(k_n)$. Suppose that the endpoint of $\tau_2$ lies
on $\alpha$. Let $\tau_2^i=M_2^i\cap \tau_2$; then, as above, 
$\alpha\cup \tau_2^i$ is a cabled graph in $M_2^i$, for $k_n$ is not
parallel to a curve lying on $S_i$.

Let $R_i$ be a regular neighborhood of $\tau_2^i$ in $M_2^i$, so that 
$R_i\cap S_i$ is a disk
$D_1^i$, and $R\cap N(k_n)$  is a disk $D_0^i$.
Assume that $D_0^i\subset A$, where $A=\eta (\alpha)$. Consider the manifold
$W_i=M_1^i\cup R_i \cup N_i(A)$, where $N_i(A)$ is a neighborhood of $A$. Let 
$\Sigma_i=\partial W_i$. As $\tau_2^i\cup \alpha$ is a cabled graph in
$M_2^i$,  it follows from [E2,6.3] that $\Sigma_i$ is incompressible in
$S^3-int\,W_i$. The  neighborhoods $R_i\cup N_i(A)$ can be chosen to be thinner
if $j>i$, that is,
$M_2^i \cap (R_j\cup N_j(A)) \subset R_i\cup N_i(A)$ if $i < j$. Then the surfaces
$\Sigma_1,\Sigma_2,\dots,\Sigma_n$ are disjoint.

Let $\tau$ be the arc obtained by sliding $\tau_1$ over $\tau_2$,
so that $M_2^i\cap \tau \subset R_i$, for all $i$.
Now take an iterate $k_{n+1}$ of $k_n$ of a special form. As before 
$k_{n+1}=\kappa_1\cup \kappa_2$, where $\kappa_2=\tau$, and $\kappa_1$ is an 
arc in $\partial N(k_n)$.
Suppose that $\kappa_1$ is contained in $A$, so that its wrapping number 
in $N_n(A)$ is $\geq 2$. Then $k_{n+1} \subset W_i$, and it 
follows from [E2,6.4] that $\Sigma_i$ is incompressible and meridionally 
incompressible in $(W_i,k_{n+1})$. So $\Sigma_i$ is a meridionally 
incompressible surface in the exterior of $k_{n+1}$ of genus $g+n+1$.
Again by [E2,8.2], it follows that $k_{n+1}$ is not parallel to a curve lying on
$\Sigma_i$

The knot $k_{n+1}$ intersects the surface $S_n$ in two points, the wrapping number
of $\kappa_2$ in $N(k_n)$ is $\geq 4$, and $k_n$ is not parallel to a curve on
$S_n$. So $k_{n+1}$ and $S_n$ satisfy the conditions of Theorem 2.1, and then
$\hat S_n=S_n\cap E(k_{n+1})$ in an incompressible, meridionally incompressible
surface in $E(k_{n+1})$ whose boundary consists of two meridians of $E(k_{n+1})$.
Then as above, by tubing $\hat S_n$ on the side of
$M_2^n$ and isotoping  into
$W_n$,  we get a closed surface $\Sigma_{n+1}$ which is incompressible but
meridionally compressible in the exterior of $k_{n+1}$. The tube added to the
surface can be chosen so that it lies in the interior of $R_n\cup N_n(A)$; this
ensures that $\Sigma_{n+1}$ is disjoint from $\Sigma_i$, for $1\leq i\leq n$.

(From the surface $S_i$ we can also get a meridionally compressible surface
$\Sigma_i^\prime$, but it will intersect $\Sigma_j$, if $i < j$. But note that
$\Sigma_{n+1}=\Sigma_n^\prime, \Sigma_{n-1}^\prime,\dots,
\Sigma_2^\prime,\Sigma_1^\prime, \Sigma_1$, are disjoint).

There is an unknotting tunnel $\beta^\prime$ for $k_{n+1}$ of the form
$\beta^\prime=k_n\cup \delta$, where $\delta$ is a straight arc in $N(k_n)$
connecting $k_n$ and $k_{n+1}$. Note that $\delta$ intersects each surface
$\Sigma_i$ in one point; this implies that 
$\Sigma_i$ is special w.r.t. $k_{n+1}$ and $\beta^\prime$. Note
also that the arc $\delta$, when going from $k_{n+1}$ to $k_n$, intersects the
surfaces in the order $\Sigma_{n+1}, \Sigma_n,\dots,\Sigma_1$.
Finally, note that the surfaces cannot be parallel, for if two of them
were, then two of the surfaces $S_i$ would also be parallel.

This shows that $k_{n+1}$, $\beta^\prime$, and
$\Sigma_{n+1},\Sigma_n,\dots,\Sigma_1$ satisfy the induction hypothesis. This
completes the proof.

By starting with a surface $S$ of genus 1, and repeating the construction 
$n-1$ times, we get the desired conclusion.
\qed
\enddemo 

\remark{Remark} It follows from the proof of the above theorem that
by changing the induction hypothesis, we can find a tunnel number one knot
$k$, with $n$ incompressible surfaces in its exterior, $S_1,S_2,\dots,S_n$, so
that $S_n$ is meridionally incompressible, but $S_i$, for $1\leq i\leq n-1$,
is meridionally compressible, and $S_n$ is the surface which is farthest
from the knot. It follows also that there are tunnel number one knots whose
exteriors contain two collection of disjoint incompressible surfaces,
$S_1,\dots,S_n$, and $\Sigma_1,\dots,\Sigma_n$, where the $S_i$ are meridionally
incompressible, and the $\Sigma_i$ are meridionally compressible. 
\endremark

\medskip
\subhead {References}\endsubhead

\enddocument